# Solving Severely Ill-Posed Linear Systems with Time Discretization Based Iterative Regularization Methods

*GONG Rongfang, HUANG Qin\**

College of Science, Nanjing University of Aeronautics and Astronautics, Nanjing 211106, P. R. China

**Abstract:** Recently, inverse problems have attracted more and more attention in computational mathematics and become increasingly important in engineering applications. After the discretization, many of inverse problems are reduced to linear systems. Due to the typical ill-posedness of inverse problems, the reduced linear systems are often ill-posed, especially when their scales are large. This brings great computational difficulty. Particularly, a small perturbation in the right side of an ill-posed linear system may cause a dramatical change in the solution. Therefore, regularization methods should be adopted for stable solutions. In this paper, a new class of accelerated iterative regularization methods is applied to solve this kind of large-scale ill-posed linear systems. An iterative scheme becomes a regularization method only when the iteration is early terminated. And a Morozov's discrepancy principle is applied for the stop criterion. Compared with the conventional Landweber iteration, the new methods have acceleration effect, and can be compared to the well-known accelerated $\nu$-method and Nesterov method. From the numerical results, it is observed that using appropriate discretization schemes, the proposed methods even have better behavior when comparing with $\nu$-method and Nesterov method.

**Key words:** linear system; ill-posedness; large-scale; iterative regularization methods; acceleration

## 0 Introduction

During the past fifty years, inverse problems have attracted more and more attention and have extensive applications in engineering and mathematical fields, such as Cauchy problem[1-3], geophysical exploration[4], steady heat conduction problems[5-6], inverse scattering problems[7-9], image processing[10-11], and reconstruction of radiated noise source[12] etc. After the discretization, many of them are reduced to the ill-posed linear system as

$$Ax = b \quad (1)$$

Here by the ill-posed linear system, it means that the condition number of the coefficient matrix $A \in \mathbf{R}^{m \times n}$ is much larger compared with the square of the scale of $A$. Denote by the quasi-solution of Eq. (1). There are three cases for $x^\dagger$: (i) The system admits a unique solution, then $x^\dagger$ is the exact solution; (ii) the system has more than one solution, then $x^\dagger$ is the one of minimal 2-norm among all solutions; (iii) the system has no solution, then $x^\dagger$ is least squares solution. It is easy to verify that $x^\dagger$ is unique.

Practically, $b$ often comes from measurements or discretization, and contains inevitably noise. Instead of $b$, assume we only have noisy data $b^\delta$ at hand, satisfying

$$\|b - b^\delta\| \leq \delta \quad (2)$$

where $\delta$ is the noisy level, and $\|\cdot\|$ is the standard Euclidean norm of a vector. Therefore, in this paper, we are devoted to find an approximate solution to the polluted linear system

$$Ax = b^\delta \quad (3)$$

Similarly, denote by $x^\delta$ the unique quasi-solution corresponding to the noisy data $b^\delta$.

For an ill-posed system, a small disturbance in $b$ will lead to much large change in solution $x$. This brings quite large difficulty when one solves the system numerically, especially when the scale is large. Thus, it is useless to use the conventional numerical methods to solve systems (1) or (3). In fact, by using the singular value decomposition (SVD), the solutions of Eq. (1) or (3) can be written formally as [13]

$$x^\dagger = \sum_{i=1}^{r} \frac{u_i^T b}{\sigma_i} v_i, \quad x^\delta = \sum_{i=1}^{r} \frac{u_i^T b^\delta}{\sigma_i} v_i$$

where $\{(\sigma_i, u_i, v_i)\}$ are the singular values and singular vectors of the coefficient matrix $A$, satisfying $\sigma_1 \geq \sigma_2 \geq \cdots \geq \sigma_r > 0, \sigma_i \leq 1$ for $i \geq k$ and some index $k \leq r$, $r$ is the rank of $A$. Therefore $\|x^\delta - x^\dagger\|$ may be quite big even if $\|b^\delta - b\|$ is small, and regularization methods are needed for obtaining a stable approximate solution to $x^\delta$ and thus also to $x^\dagger$.

Generally speaking, there are three groups of regularization methods: truncated singular value decomposition (TSVD)[13-15], variational regularization methods[16-17] and iterative regularization methods[18-21]. For large-scale ill-posed problems, TSVD leads to computational difficulty in that the vast storage and the heavy computing burden. The most famous variational regularization method is Tikhonov's regularization method [22]. Recently, a projection fractional Tikhonov regularization method is proposed[23]. However, during the determination of the regularization parameters, a forward problem has to be solved for each regularization parameter, which makes the calculation very large. Iterative regularization methods have the advantages of low computational cost and simple forms. Thus for large-scale problem, we prefer to use iterative regularization methods, where iterative schemes are proposed to solve the following optimization problem.

$$\min J(x) = \frac{1}{2} \|Ax - b^\delta\|^2 \tag{4}$$

The most classic iterative regularization method should be Landweber iteration[24-25]. For the linear system (3), Landweber iteration is defined by

$$x_{k+1} = x_k - \Delta t \nabla J(x_k) \tag{5}$$

where $\nabla J(x_k) = A^T(Ax_k - b^\delta)$ and $\Delta t \in (0, 2/\|A\|^2)$.

Eq. (5) can be viewed as a discrete analog of the following first order evolution equation.

$$\begin{cases} \dot{x}(t) = -\nabla J(x(t)) \\ x(0) = x_0 \end{cases} \tag{6}$$

where $t$ is the introduced artificial time, and $x_0$ is an initial guess to $x^\dagger$. The formulation (6) is known as an asymptotical regularization, or the Showalter's method[26-27]. It is well known that the Landweber method converges quite slowly[28]. Thus, it is no wonder that accelerating strategies are adopted in practice. In recent years, there has been increasing evidence to show that the second order iterative methods exhibit remarkable acceleration properties for ill-posed problems. The most well-known methods are the Levenberg–Marquardt method[29], the iteratively regularized Gauss–Newton method[30], the Nesterov acceleration scheme[31] and the $\nu$-method[28, § 6.3]. Here, the iterative schemes are proposed by discretizing the following second order evolution system[32].

$$\begin{cases} \ddot{x}(t) + \eta(t)\dot{x}(t) + A^T A x(t) = A^T b^\delta \quad t \in (t_0, \infty) \\ x(t_0) = x_0, \dot{x}(t_0) = \dot{x}_0 \end{cases} \tag{7}$$

where $(x_0, \dot{x}_0)$ is the prescribed initial data, $\eta > 0$, uniformly bounded, is the damping parameter, and may or may not depends on the time $t$. Note that the system (7) becomes a regularization method only when the evolution process is stopped in advance. Denote by $T$ the terminating time. Then it plays the role of the regularization parameter and should be chosen properly. Here, the Morozov's discrepancy principle is as the stopping criterion

$$\|Ax(T) - b^\delta\| \leq \tau\delta \tag{8}$$

where $\tau$ is a fixed positive number. Like Landweber method, the well-known Nesterov acceleration method can be viewed as the discrete analog of Eq. (7). In fact, for all fixed $T > 0$ [27]

$$\lim_{\omega \to 0} \max_{0 \leq k \leq T/\sqrt{\omega}} \|x_k - x(k\sqrt{\omega})\| = 0$$

where $x(\cdot)$ is the solution of Eq. (7) with $\eta(t) = \alpha/t$, and $\{x_k\}$ is the sequence, generated by the Nesterov's scheme with parameters $(\alpha, \omega)$.

The purpose of this paper is to apply the second order dynamical system (7), together with the discrepancy principle (8) for the choice of the terminating time, to the ill-posed linear system (3). Three aspects are expected to be addressed: (i) On one hand, unlike the existing work, where the exact solution is often assumed to exist, the system here may have no solution, and thus the existing theoretical results could not be used directly; on the other hand, benefiting by the linearity and the finite dimension of the problem, compared with the theoretical analysis in the literature, arguments here can be largely simplified. (ii) Effect of the magnitude of the problem on the iterations of different methods is investigated. (iii) Effect of the ill-posedness extent of the problem on the iterations of different methods is discussed.

# 1 Theoretical Analysis of Continuous Second-Order Flow

In this section, we are devoted to give a series of theoretical analysis on the second order dynamical system (7). Without loss of generality, we set $t_0 = 0$ when $\eta =$ constant, and $t_0 = 1$ when $\eta$ is time-dependent; set $\dot{x}_0 = 0$ when $\eta$ is time-dependent. Moreover, for a dynamical damping parameter, we take $\eta=(1+2s)/t$ (the constant $s > -1/2$) as an example for the theoretical analysis. Theoretical results for other choice such as $\eta = 1/\sqrt{t}$ could be analyzed similarly, and we omit here.

The definition of a regularized solution is first introduced.

**Definition** Let $x(t) \in \mathrm{R}^n$ be the solution of System (7). Then, $x(T^\delta)$, equipped with an appropriate terminating time $T^\delta = T(\delta, b^\delta)$, is called a second order asymptotical regularization solution of Eq. (7) if $x(T^\delta)$ converges (strongly) to $x^\dagger$ in $\mathrm{R}^n$ as $\delta \to 0$.

## 1.1 Existence and uniqueness of solution trajectory

About the global existence and uniqueness of solution to Eq. (7), the following result can be obtained.

**Theorem 1** For any pair $(x_0, \dot{x}_0) \in \mathrm{R}^n \times \mathrm{R}^n$, there exists a unique solution $x(\cdot) \in C^\infty([t_0, \infty), \mathrm{R}^n)$ for the second order dynamical system (7). Moreover, $x$ depends continuously on the data $b^\delta$.

**Proof** Denote $z = (x, \dot{x})^\mathrm{T}, z_0 = (x_0, \dot{x}_0)^\mathrm{T}$ and rewrite Eq. (7) as a first order differential system.

$$\dot{z}(t) = Bz(t) + \tilde{b} \triangleq F(z(t)) \qquad (9)$$

where $B = \begin{pmatrix} 0 & I \\ -A^\mathrm{T}A & -\eta(t)I \end{pmatrix}$, $\tilde{b} = \begin{pmatrix} 0 \\ A^\mathrm{T}b^\delta \end{pmatrix}$ and $I$ denotes the identity matrix of order $n$. Since $A$ is a bounded matrix and $\eta$ is uniformly bounded, $B$ is a uniformly bounded matrix. Hence, $F(\cdot)$ is a global Lipschitz functional. By the Cauchy-Lipschitz-Picard theorem [33-34], the first order autonomous system (9) has a global unique solution $z$ for any given initial data $z_0$, and then the second order dynamical system (7) has a global unique solution $x$ for any given initial data $(x_0, \dot{x}_0)^\mathrm{T}$. Since $A^\mathrm{T}(b^\delta - Ax)$ is linear with respect to $x$, and $\eta = $ constant or $(1+2s)/t$, both of them are infinitely differentiable, which gives $x(\cdot) \in C^\infty([t_0, \infty), \mathrm{R}^n)$. Moreover, the continuous dependence of $x$ on the data $b^\delta$ can be easily verified. ∎

Next the relationship between the solution $x(t)$ of Eq. (7) and the exact one $x^\dagger$ is shown. For the future use, denote by $S$ the set of minimizers of $J(\cdot)$ which can be characterized by

$$S = \{x \in \mathrm{R}^n \mid \nabla J(x) = 0\}$$

It is easy to conclude that $S$ is a non-void closed convex subset of $\mathrm{R}^n$. Moreover, for the statement of clarity, the discussion is divided into two parts, corresponding to the noise-free data $b$ and noisy data $b^\delta$ respectively.

## 1.2 Limiting behavior of the solution for noise-free data

Define the modified energy functional of $x$ as

$$E(t) = J(x(t)) - J(x^\dagger) + \frac{1}{2}\|\dot{x}(t)\|^2 \qquad (10)$$

**Case I** $\eta = $ constant and $t_0 = 0$, that is, we consider the following system

$$\begin{cases} \ddot{x}(t) + \eta\dot{x}(t) + A^\mathrm{T}Ax(t) = A^\mathrm{T}b & t \in (0, \infty) \\ x(0) = x_0, \dot{x}(0) = \dot{x}_0 \end{cases} \qquad (11)$$

**Proposition 1** Let $x$ be the solution of Eq. (11). Then, the following properties hold

(i) $x(\cdot) \in L^\infty([0, \infty), \mathrm{R}^n)$

(ii) $\dot{x}(\cdot) \in L^\infty([0, \infty), \mathrm{R}^n) \cap L^2([0, \infty), \mathrm{R}^n)$ and thus $\dot{x}(t) \to 0$ as $t \to \infty$

(iii) $\ddot{x}(\cdot) \in L^\infty([0, \infty), \mathrm{R}^n)$ and $\ddot{x}(t) \to 0$ as $t \to \infty$

(iv) $\lim_{t \to \infty} \nabla J(x(t)) = 0$.

(v) $J(x(t)) - J(x^\dagger) = o(t^{-1})$ as $t \to \infty$

**Proof** (i) Differentiating $E(t)$ of Eq. (10) and using the system (11) to give

$$\dot{E}(t) = -\eta\|\dot{x}(t)\|^2 \qquad (12)$$

which means $E(t)$ is non-increasing, and thus
$$E(t) \leq E(0) \tag{13}$$
holds for all $t \geq 0$. Consequently, $x(\cdot) \in L^\infty([0,\infty), \mathrm{R}^n)$ by combining Eq. (13) and the coerciveness of $J(\cdot)$, i.e.
$$\lim_{\|x\| \to +\infty} J(x) = +\infty$$

(ii) On one hand, due to Eq. (13), we have
$$\|\dot{x}(t)\|^2 \leq 2(E(t) + J(x^\dagger)) \leq 2(E(0) + J(x^\dagger))$$
which gives $\dot{x}(\cdot) \in L^\infty([0,\infty), \mathrm{R}^n)$. On the other hand, by using the decrease and non-negativity of $E(t)$, the limit $E_\infty := \lim_{t \to \infty} E(t)$ exists. Therefore, integrating both sides in Eq. (12), we obtain
$$\int_0^\infty \|\dot{x}(t)\|^2 dt = \frac{E(0) - E_\infty}{\eta} < \infty$$
which yields $\dot{x}(\cdot) \in L^2([0,\infty), \mathrm{R}^n)$.

In addition, according to a classical result, $\dot{x}(\cdot) \in L^\infty([0,\infty), \mathrm{R}^n) \cap L^2([0,\infty), \mathrm{R}^n)$ implies $\dot{x}(t) \to 0$ as $t \to \infty$.

(iii) From Eq. (11), we obtain
$$\ddot{x}(t) = -\eta \dot{x}(t) + A^\mathrm{T}(b - Ax(t))$$
which gives immediately that $\ddot{x}(\cdot) \in L^\infty([0,\infty), \mathrm{R}^n)$.

By differentiating the first equation of Eq. (11), we obtain
$$\dddot{x}(t) + \eta \ddot{x}(t) = -A^\mathrm{T} A \dot{x}(t) \triangleq g(t) \tag{14}$$

Denote by $y = \ddot{x}$. Then Eq. (14) is reduced to
$$\dot{y}(t) + \eta y(t) = g(t)$$
which implies $y \to 0$ as $t \to \infty$ by noticing that $\eta > 0$ and $g(t) \to 0$ as $t \to \infty$. Thus $\ddot{x}(t) \to 0$ as $t \to \infty$.

(iv) By using Eq. (11) again, we have
$$\nabla J(x(t)) = A^\mathrm{T}(Ax(t) - b) = -\ddot{x}(t) - \eta \dot{x}(t)$$
which gives $\lim_{t \to \infty} \nabla J(x(t)) = 0$ by using the facts in properties (ii) and (iii).

(v) Define
$$h(t) = \frac{\eta}{2} \|x(t) - x^\dagger\|^2 + \langle \dot{x}(t), x(t) - x^\dagger \rangle \tag{15}$$

By elementary calculations, we derive that
$$\dot{h}(t) = \eta \langle \dot{x}(t), x(t) - x^\dagger \rangle + \langle \ddot{x}(t), x(t) - x^\dagger \rangle + \|\dot{x}(t)\|^2$$
$$= \|\dot{x}(t)\|^2 - \|Ax(t) - b\|^2$$

which implies that (by noting $\dot{E}(t) = -\eta \|\dot{x}(t)\|^2$)
$$\dot{E}(t) + \eta E(t) + \frac{\eta}{2} \dot{h}(t) = 0$$
or
$$\eta E(t) = -\dot{E}(t) - \frac{\eta}{2} \dot{h}(t)$$

Integrate the above equation on $[0,T]$ to obtain, together with the non-negativity of $E(t)$,
$$\begin{aligned} \int_0^T E(t) dt &= \frac{1}{\eta}(E(0) - E(T)) + \frac{1}{2}(h(0) - h(T)) \\ &\leq \frac{1}{\eta} E(0) + \frac{1}{2}(h(0) - h(T)) \end{aligned} \tag{16}$$

Due to propertes (i)-(ii), $x(t), \dot{x}(t)$ are uniformly bounded, and so is $h(t)$. Hence, letting $T \to \infty$ in Eq. (16), we obtain
$$\int_0^\infty E(t) dt < \infty \tag{17}$$

Since $E(t)$ is non-increasing, we deduce that
$$\int_{T/2}^T E(t) dt \geq \frac{T}{2} E(T) \geq \frac{T}{2} \left( J(x(T)) - J(x^\dagger) \right). \tag{18}$$

Using Eq. (17), the left side of Eq. (18) tends to 0 when $T \to \infty$, which implies that $\lim_{t \to \infty} t \left( J(x(t)) - J(x^\dagger) \right)/2 = 0$ or $J(x(t)) - J(x^\dagger) = o(t^{-1})$ as $t \to \infty$. The proof is completed. ∎

**Case II** $\eta(t) = (1+2s)/t$ and $t_0 = 1$

Now, we consider the following evolution system
$$\begin{cases} \ddot{x}(t) + \dfrac{1+2s}{t} \dot{x}(t) + A^\mathrm{T} A x(t) = A^\mathrm{T} b \quad t \in (1, \infty) \\ x(1) = x_0, \dot{x}(1) = 0 \end{cases} \tag{19}$$

**Proposition 2** Let $x$ be the solution of Eq. (19). Then, for $s \geq 1$, the following properties hold
(i) $x(\cdot) \in L^\infty([1,\infty), \mathrm{R}^n)$
(ii) $\dot{x}(\cdot) \in L^\infty([1,\infty), \mathrm{R}^n) \cap L^2([1,\infty), \mathrm{R}^n)$ and thus
$$\dot{x}(t) \to 0 \text{ as } t \to \infty$$
(iii) $\ddot{x}(\cdot) \in L^\infty([1,\infty), \mathrm{R}^n)$
(iv) $J(x(t)) - J(x^\dagger) = \mathcal{O}(t^{-2})$ as $t \to \infty$

**Proof** (i) Differentiating $E(t)$ of Eq. (10) and using the system (19) to give
$$\dot{E}(t) = -\frac{1+2s}{t} \|\dot{x}(t)\|^2 \tag{20}$$

The remaining proof for property (i) is similar to that for property (i) of proposition 1.

(ii) Like the statements in the proof of proposition 1,

we have $\dot{x}(\cdot) \in L^\infty([1,\infty), \mathbb{R}^n)$, the limit $E_\infty := \lim_{t\to\infty} E(t)$ exists. Now we show $\dot{x}(\cdot) \in L^2([1,\infty), \mathbb{R}^n)$. To the end, define

$$E_1(t) = t^2(J(x(t)) - J(x^\dagger)) + \frac{1}{2}\left\|2(x(t) - x^\dagger) + t\dot{x}(t)\right\|^2 + 2(s-1)\left\|x(t) - x^\dagger\right\|^2$$

By using the convexity inequality $J(x^\dagger) \geq J(x) + (\nabla J(x), x^\dagger - x)$ for all $x \in \mathbb{R}^n$ and Eq. (19), it is not difficult to show that

$$\dot{E}_1(t) \leq -2(s-1)t\|\dot{x}(t)\|^2 \quad (21)$$

Hence, for $s \geq 1$, $\dot{E}_1(t) \leq 0$ and thus $E_1(\cdot)$ is non-increasing. Together with the fact that $E_1(t) \geq 0$ for all $t \geq 0$, the limit $E_1(\infty) := \lim_{t\to\infty} E_1(t)$ exists.

Integrating both sides in Eq. (21), we obtain that

$$\int_1^\infty \|\dot{x}(t)\|^2 dt \leq \int_1^\infty t\|\dot{x}(t)\|^2 dt \leq \frac{E_1(1) - E_1(\infty)}{2(s-1)} < \infty \quad (22)$$

which yields $\dot{x}(\cdot) \in L^2([1,\infty), \mathbb{R}^n)$. Combining it and $\dot{x}(\cdot) \in L^\infty([1,\infty), \mathbb{R}^n)$ to conclude that $\dot{x}(t) \to \mathbf{0}$ as $t \to \infty$.

(iii) From Eq. (19), we obtain

$$\ddot{x}(t) = -\frac{1+2s}{t}\dot{x}(t) + A^T(b - Ax(t))$$

which gives immediately that $\ddot{x}(\cdot) \in L^\infty([1,\infty), \mathbb{R}^n)$.

(iv) Define

$$E_2(t) = t^2(J(x(t)) - J(x^\dagger)) + \frac{1}{2}\|2s(x(t) - x^\dagger) + t\dot{x}(t)\|^2 \quad (23)$$

By using the convexity inequality $J(x^\dagger) \geq J(x) + (\nabla J(x), x^\dagger - x)$ and Eq. (19) again, we have

$$\dot{E}_2(t) \leq -2(s-1)t\left(J(x(t)) - J(x^\dagger)\right)$$

Similar to $E_1(\cdot)$, for $s \geq 1$, $E_2(\cdot)$ is non-increasing, and the limit $\lim_{t\to\infty} E_2(t)$ exists. Therefore, from Eq. (23), there holds

$$t^2\left(J(x(t)) - J(x^\dagger)\right) \leq E_2(t) \leq E_2(1)$$

or

$$0 \leq J(x(t)) - J(x^\dagger) \leq E_2(1)t^{-2}$$

which gives $J(x(t)) - J(x^\dagger) = \mathcal{O}(t^{-2})$ as $t \to \infty$. The proof is completed. ∎

**Remark 1** Compared with the first order method (6), where the convergence rate of the objective functional is $J(x(t)) - J(x^\dagger) = \mathcal{O}(t^{-1})$, the convergence rates $J(x(t)) - J(x^\dagger) = o(t^{-1})$ in proposition 1 and $J(x(t)) - J(x^\dagger) = \mathcal{O}(t^{-2})$ in proposition 2 show that the second order dynamical systems (11) and (19) can achieve higher convergence order, indicating the second order dynamical system has a property of acceleration.

Now we are in a position to give the convergence result of the dynamical solution $x(t)$ of Eq. (7) with the noise-free data. Recall that $x^\dagger$ is the minimum norm minimizer of Eq. (4) with the exact data $b$. Then, formally, we have

$$x^\dagger = (A^T A)^+ A^T b$$

where the superscript + means the Moore-Penrose inverse.

**Theorem 2** Let $x(\cdot)$ be the solution of Eq. (7) with the exact data, then

$$\lim_{t\to\infty} x(t) = x^\dagger. \quad (24)$$

**Proof** Let $\{\sigma_j; u_j, v_j\}_{j=1}^\infty$ be the singular system for matrix $A$, i.e. we have $Av_j = \sigma_j u_j$ and $A^T u_j = \sigma_j v_j$ with ordered singular values $\|A\| = \sigma_1 \geq \sigma_2 \geq \cdots \geq \sigma_r > 0$.

Let's first consider the case with fixed damping parameter $\eta$. Similar to Ref. [32], we distinguish three different cases: (a) The overdamped case: $\eta > 2\|A\|$, (b) the underdamped case: there is an index $j_0$ such that $2\sigma_{j_0+1} < \eta < 2\sigma_{j_0}$, and (c) the critical damping case: an index $j_0$ exists such that $\eta = 2\sigma_{j_0}$. For simplicity, we only consider the overdamped case. The other two cases can be studied similarly. According to Ref. [32], we have

$$x(t) = (I - A^T A g(t, A^T A))x_0 + \phi(t, A^T A)\dot{x}_0 + g(t, A^T A)A^T b \quad (25)$$

where

$$\begin{cases} g(t,\lambda) = \frac{1}{\lambda}\left(1 - \frac{\lambda_1}{\sqrt{\eta^2 - 4\lambda}}e^{-\lambda_2 t} + \frac{\lambda_2}{\sqrt{\eta^2 - 4\lambda}}e^{-\lambda_1 t}\right) \\ \phi(t,\lambda) = -\frac{1}{2\sqrt{\eta^2 - 4\lambda}}\left(e^{-\lambda_2 t} - e^{-\lambda_1 t}\right) \end{cases} \quad (26)$$

and $\lambda_1 = \frac{\eta + \sqrt{\eta^2 - 4\lambda}}{2}, \lambda_2 = \frac{\eta - \sqrt{\eta^2 - 4\lambda}}{2}$.

Together with theorem 1, it is straightforward to check that $x(t)$ defined Eq. (25) is a unique solution of Eq. (11). Moreover, we can deduce from Eq. (25) that $x(t) \to (A^T A)^+ A^T b$ as $t \to \infty$ by noting that

$$1 - \lambda g(t,\lambda) = \frac{\lambda_1}{\sqrt{\eta^2 - 4\lambda}} e^{-\lambda_2 t} - \frac{\lambda_2}{\sqrt{\eta^2 - 4\lambda}} e^{\lambda_1 t}$$

Next we consider the case with dynamic damping parameter $\eta(t) = (1+2s)/t$. According to Ref. [35], we have (note that $\dot{x}_0 = 0$)

$$x(t) = (I - A^T A g_2(t, A^T A))x_0 + g_2(t, A^T A)A^T b, \quad (27)$$

Where

$$g_2(t,\lambda) = \frac{1 - 2^s \Gamma(s+1) \frac{J_s(\sqrt{\lambda}t)}{(\sqrt{\lambda}t)^s}}{\lambda}, \quad (28)$$

$\Gamma(\cdot)$ and $J_s(\cdot)$ denote the Gamma function and the Bessel function of first kind of order $s$, respectively. As in the first case, $x(t)$ defined in Eq. (27) is the unique solution of Eq. (19). By using the asymptotic

$$J_s(\sqrt{\lambda}t) = \sqrt{\frac{2}{\pi}}(\sqrt{\lambda}t)^{-\frac{1}{2}} \cos(\sqrt{\lambda}t - \frac{\pi(2s+1)}{4}) + \mathcal{O}(\frac{1}{\sqrt{\lambda}t}) \quad t \to \infty \quad (29)$$

for any fixed $\lambda > 0$, we conclude that $x(t) \to (A^T A)^+ A^T b$ as $t \to \infty$. The proof is completed. ∎

## 1.3 Regularization properties for the noisy data

In section 1.2, we study the limiting behaviors of the solution for the second dynamical system (7) with noise-free data. It is shown that for exact data $b$, theoretically, the larger the time $t$ is, the more accurate the solution $x(t)$ is. Therefore, we are readily given a satisfactory approximate solution $x(T^*)$ to $x^\dagger$ for a large enough time $T^*$. However, due to the inevitable noise, practically, it is not the case. In the case that only a noisy data $b^\delta$ is available, the system (7) cannot produce a reasonable approximate solution unless the terminating time $T^*$ is chosen appropriately. In other word, the asymptotic process should be terminated in advance.

In this subsection, we investigate the regularization property of the dynamic solution $x^\delta(t)$ of system (7), equipped with the Morozov's discrepancy principle for selecting the terminating time $T^*$.

We define the Morozov's discrepancy function as

$$\chi(T) = \| Ax^\delta(T) - b^\delta \| - \tau\delta \quad (30)$$

where $x^\delta(t)$ is the solution of system (7), and $\tau$ is a fixed positive number.

**Proposition 3** For any $\tau > 0$, let the initial guess $x_0$ is chosen satisfying $\|Ax_0 - b^\delta\| \geq \tau\delta$. Then $\chi(T)$ has at least one zero point.

**Proof** First, by theorem 1, it is easy to conclude that $\chi(\cdot)$ is continuous over $[t_0, \infty)$. Recall that $x^\delta$ is the minimum norm minimizer of Eq. (4). Then, formally, like $x^\dagger$ we have

$$x^\delta = (A^T A)^+ A^T b^\delta$$

Applying theorem 2 to $x^\delta(t)$ yield

$$\lim_{t \to \infty} x^\delta(t) = x^\delta$$

which implies

$$\lim_{T \to \infty} \chi(T) = \|Ax^\delta - b^\delta\| - \tau\delta = -\tau\delta < 0 \quad (31)$$

By combing Eq. (31) and $\chi(t_0) = \|Ax_0 - b^\delta\| - \tau\delta > 0$ as well as the continuity, we conclude that the nonlinear equation $\chi(T) = 0$ admits at least one root.

Now we are in a position to present a convergence result in the following.

**Theorem 3** Let $x^\delta(t)$ be the solution of the dynamic system (7). Denote by $T^* = T(\delta, b^\delta)$ the first zero point of $\chi(T)$. Then $x^\delta(T^*)$ converges to $x^\dagger$ in $\mathbb{R}^n$ as $\delta \to 0$.

**Proof** As $\delta \to 0$, $T(\delta, b^\delta)$ may go to infinity or may have a finite accumulation point. We first consider the situation that $T(\delta, b^\delta) \to \infty$. Recall that $A$ is a $m \times n$ matrix with rank $r$. We proceed to the proof with two possible cases: (i) $r \geq n$; (ii) $r < n$.

(i) In the case $r \geq n$, $A^T A$ is full rank, and thus

$$x^\dagger = (A^T A)^{-1} A^T b, \quad x^\delta = (A^T A)^{-1} A^T b^\delta \quad (32)$$

By using triangle inequality,

$$\|x^\delta(T^*) - x^\dagger\| \leq \|x^\delta(T^*) - x^\delta\| + \|x^\delta - x^\dagger\| \quad (33)$$

On one hand, by applying the convergence in theorem 2, we have

$$\lim_{\delta \to 0} \|x^\delta(T^*) - x^\delta\| = 0 \quad (34)$$

On the other hand, from Eq. (32),

$$\|x^\delta - x^\dagger\| \leq \|(AA^T)^{-1}\| \|A^T(b^\delta - b)\| \leq \frac{1}{\|A\|}\delta \quad (35)$$

Combine Eq. (33)-(35) to conclude that
$$\lim_{\delta \to 0} \left\| x^{\delta}(T^{*}) - x^{\dagger} \right\| = 0$$

(ii) For the case $r < n$, both the systems (1) and (3) have infinite solution, and $x^{\dagger}, x^{\delta}$ are the minimum norm solutions, respectively. Let $\{\delta_n\}$ be a sequence converging to 0 as $n \to \infty$, and let $b^{\delta_n}$ be a corresponding sequence of noisy data with $\left\| b^{\delta_n} - b \right\| \leq \delta_n$. For a pair $(\delta_n, b^{\delta_n})$, denote by $T_n^*$ the corresponding terminating time point determined from the discrepancy principles $\chi(T) = 0$. According to the continuity of $x^{\delta_n}(t)$, for any $\varepsilon > 0$ and large enough $n$ and thus $T_n^*$, there exists a point $T < T_n^*$ such that
$$\left\| x^{\delta_n}(T_n^*) - x^{\delta_n}(T) \right\| \leq \varepsilon/3 \tag{36}$$

From the continuity of $x^{\delta_n}(\cdot)$ due to theorem 2, for large enough $n$, there holds
$$\left\| x^{\delta_n}(T) - x(T) \right\| \leq \varepsilon/3 \tag{37}$$

Moreover, by applying theorem 2 again, we obtain, for large enough $T$ (because $T_n^*$ could be large enough),
$$\left\| x(T) - x^{\dagger} \right\| \leq \varepsilon/3 \tag{38}$$

Then by combing Eqs. (36)-(38), and using triangle inequality
$$\left\| x^{\delta_n}(T_n^*) - x^{\dagger} \right\| \leq \left\| x^{\delta_n}(T_n^*) - x^{\delta_n}(T) \right\|$$
$$+ \left\| x^{\delta_n}(T) - x(T) \right\| + \left\| x(T) - x^{\dagger} \right\|$$

There is a number $n_1$ such that for any $n \geq n_1$
$$\left\| x^{\delta_n}(T_n^*) - x^{\dagger} \right\| \leq \varepsilon.$$

Since $\varepsilon$ is arbitrary, we arrive at the convergence of $x^{\delta_n}(T_n^*)$ to $x^{\dagger}$ as $n \to \infty$.

For the case that $T(\delta, b^{\delta})$ has a finite accumulation point, we can use arguments similar to those in Theorem 2.4 of Ref. [36], and omit here. The proof is completed. ∎

## 2 Iterative Schemes

For the numerical implementation, this section is devoted to present several iterative schemes for the resolution of problem (4). To the end, we first convert the system (7) into
$$\begin{cases} \dot{x}(t) = q(t), \\ \dot{q}(t) = -\eta(t) q(t) + A^{\mathrm{T}}(b^{\delta} - A x(t)), \\ x(t_0) = x_0, \dot{x}(t_0) = \dot{x}_0, \end{cases} \tag{39}$$

or

$$\begin{cases} \dot{y}(t) \\ = \begin{pmatrix} 0 & I \\ -A^{\mathrm{T}} A & -\eta(t) I \end{pmatrix} \begin{pmatrix} x(t) \\ q(t) \end{pmatrix} + \begin{pmatrix} 0 \\ A^{\mathrm{T}} b^{\delta} \end{pmatrix} \\ \triangleq A(t) y(t) + \tilde{b} \\ x(t_0) = x_0, \dot{x}(t_0) = \dot{x}_0 \end{cases}$$
(40)

Then, the numerical discretization of the differential system (39) or (40) together with the discrepancy principle produces second order iterative regularization methods. The damped symplectic integrators are extremely attractive since these schemes are closely related to the canonical transformations[37], and the trajectories of the discretized second flow are usually more stable for its long-term performance. Applying the Symplectic Euler method to the system (39), we obtain the iterative scheme at the $k$ th step
$$\begin{cases} q^{k+1} = q^k + \Delta t \left( A^{\mathrm{T}}(b^{\delta} - A x^k) - \eta(t_k) q^k \right) \\ x^{k+1} = x^k + \Delta t q^{k+1} \\ x^0 = x_0, q^0 = \dot{x}_0 \end{cases} \tag{41}$$

where $x^k = x^{\delta}(t_k), q^k = q(t_k)$, and $\Delta t$ is the step-size. By elementary calculations, Scheme (41) can be expressed as the form of following three-term semi-iterative method
$$x^{k+1} = x^k + a_k (x^k - x^{k-1}) + \omega_k A^{\mathrm{T}}(b^{\delta} - A x^k) \tag{42}$$
where $a_k = (1 - \Delta t \eta(t_k))$ and $\omega_k = (\Delta t)^2$.

For the high order Symplectic methods, the Störmer-Verlet scheme can be considered as
$$\begin{cases} q^{k+1/2} = q^k - \dfrac{\Delta t}{2} \eta(t_k) q^{k+1/2} + \dfrac{\Delta t}{2} A^{\mathrm{T}}(b^{\delta} - A x^k), \\ x^{k+1} = x^k + \Delta t q^{k+1/2}, \\ q^{k+1} = q^{k+1/2} - \dfrac{\Delta t}{2} \eta(t_{k+1}) q^{k+1/2} + \dfrac{\Delta t}{2} A^{\mathrm{T}}(b^{\delta} - A x^{k+1}), \\ x^0 = x_0, q^0 = \dot{x}_0. \end{cases} \tag{43}$$

Like Eq. (41), the iteration (43) can also be rewritten in the form of Eq. (42), but with parameters
$$a_k = \dfrac{1 - \dfrac{\Delta t}{2} \eta(t_k)}{1 + \dfrac{\Delta t}{2} \eta(t_k)}, \omega_k = \dfrac{(\Delta t)^2}{1 + \dfrac{\Delta t}{2} \eta(t_k)}$$

In Ref. [35], a modified Störmer-Verlet scheme is proposed as

$$\begin{cases} q^{k+1/2} = q^k - \frac{\Delta t}{2}\eta(t_k)q^{k+1/2} + \frac{\Delta t}{2}A^T(b^\delta - Ax^k) \\ x^{k+1} = x^k + \Delta t q^{k+1/2} \\ v^{k+1} = x^{k+1} + 2\Delta t a_{k+1} q^{k+1/2} \\ q^{k+1} = q^{k+1/2} - \frac{\Delta t}{2}\eta(t_{k+1})q^{k+1/2} + \frac{\Delta t}{2}A^T(b^\delta - Av^{k+1}) \\ x^0 = x_0, q^0 = \dot{x}_0 \end{cases} \quad (44)$$

with $a_k = \frac{1-\Delta t \eta(t_k)}{1+\Delta t n(t_k)}$. The third step in Eq. (44) is inspired by the Nesterov's method. Unlike the Symplectic Euler method and Störmer-Verlet scheme, the scheme (44) expresses the recurrence form as

$$x^{k+1} = x^k + a_k(x^k - x^{k-1}) \\ + \omega_k A^T(b^\delta - A(x^k + a_k(x^k - x^{k-1}))), k=1,2,\cdots, \quad (45)$$

with parameters

$$a_k = \frac{1-\frac{\Delta t}{2}n(t_k)}{1+\frac{\Delta t}{2}\eta(t_k)}, \omega_k = \frac{(\Delta t)^2}{1+\frac{\Delta t}{2}\eta(t_k)} \quad (46)$$

As indicated by theorems 2 and 3, theoretically, for the noise-free data or noisy data with small noisy level, the iteration should go far enough for a good enough approximate solution $x^k$. In this case, it requires long time behavior for the iteration. Keep in mind that we are solving ill-posed problems. Therefore, when the data contains noise with not small noisy level, the iteration should stop in advance before the solution gets worse. Therefore, the accuracy itself in the discretization of system (39) or (40) plays more important role than the requirement of the long behavior. We can apply high order Runge-Kutta methods for this purpose. As an example, the classical fourth order Runge-Kutta method is adopted as

$$\begin{cases} K_1 = A(t_k)y_k + \tilde{b}, \\ K_2 = A(t_{k+1/2})(y_k + \Delta t K_1/2) + \tilde{b}, \\ K_3 = A(t_{k+1/2})(y_k + \Delta t K_2/2) + \tilde{b}, \\ K_4 = A(t_{k+1})(y_k + \Delta t K_3) + \tilde{b}, \\ y_{k+1} = y_k + \frac{\Delta t}{6}(K_1 + 2K_2 + 2K_3 + K_4), \\ y_0 = (x_0, \dot{x}_0)^T. \end{cases} \quad (47)$$

We note that other Runge-Kutta methods are possible.

After the discretization, the regularization parameter is reduced to the iterative step $k$. For the schemes (41), (43), (44) and (47), the Morozov's discrepancy principle reads: find the smallest $k$ satisfying

$$\|Ax_k - b^\delta\| \leq \tau\delta \quad (48)$$

## 3 Numerical Results

Two numerical examples are presented to demonstrate the effectiveness of the iterative regularization methods (41), (43), (44) and (47) for the resolution of the system (3). For the sake of simplicity, the Symplectic Euler method (41) is termed as SE1 or SE2, corresponding to a constant or time-dependent damping parameter respectively. Similarly, the Störmer-Verlet method is termed as SV1 or SV2, the Modified Störmer-Verlet method is termed as MSV1 or MSV2, and the Runge-Kutta method is termed as RK1 or RK2. In the following experiments, the maximal iterative number $N_{max}$ is set as 5000.

For the comparison of the proposed methods with the existing work, we introduce such classical iteration methods as Landweber method, conjugate gradient method (termed as CG), $\nu$-method and the Nesterov's method. For the problem (3), the Landweber method (5) reads

$$x_{k+1} = x_k - \Delta t A^T(Ax_k - b^\delta) \quad (49)$$

with $0 < \Delta t < 2/\|A\|^2$. The $\nu$-method has the form

$$x^{k+1} = x^k + \mu_{k+1}(x^k - x^{k-1}) - \omega_{norm}\omega_k g_k \quad (50)$$

with $\mu_1 = 0, \omega_1 = (4\nu+2)/(4\nu+1)$, $k=1,2,\cdots$

$$\mu_k = \frac{(k-1)(2k-3)(2k+2\nu-1)}{(k+2\nu-1)(2k+4\nu-1)(2k+2\nu-3)},$$

$$\omega_k = 4\frac{(2k+2\nu-1)(k+\nu-1)}{(k+2\nu-1)(2k+4\nu-1)},$$

where $x^{-1}, x^0$ are the initial values, and $g_k = A^T(Ax^k - b^\delta)$. In (50), $\omega_{norm}$ plays the role of normalization, and it can be set as $\omega_{norm} = 1/\|A\|^2$.

The Nesterov's method is defined by

$$x^{k+1} = x^k + \frac{k-1}{k+\alpha-1}(x^k - x^{k-1}) - \omega g_k \quad (51)$$

with $\alpha > 3$, $0 < \omega \leq \omega_{norm}$ and $k=0,1,2,\cdots$.

Like the methods in section 2, the Morozov's discrepancy principle (48) is applied to Landweber method, CG method, $\nu$-method and the Nesterov method for the choice of iterative numbers.

### 3.1. Example 1

In the first example, a Fredholm integral equation of

the first kind of convolution type in one space dimension is considered as[38-39]

$$(\mathcal{K}f)(x) = \int_0^1 k(x-x')f(x')d(x') = g(x) \quad (52)$$
$$0 < x < 1$$

where the kernel function $k(x) = C\exp(-x^2/2\gamma^2)$ with positive constants $C$ and $\gamma$. Using numerical integration, Eq. (52) is discretized to a linear system

$$\boldsymbol{Kf} = \boldsymbol{d} \quad (53)$$

where $[\boldsymbol{K}]_{ij} = hC\exp\left(-\dfrac{((i-j)h)^2}{2\gamma^2}\right)$, $[\boldsymbol{f}]_i = f(ih)$,

$[\boldsymbol{d}]_i = g(ih)$, $h = 1/n$, $1 \le i,j \le n$. In the experiments of this subsection, fix $\gamma=0.05$, $C=1/\gamma$. The dependence of the condition number of the matrix $\boldsymbol{K}$ on the magnitude is plotted in Fig. 1 which shows that $\boldsymbol{K}$ is ill-posed as long as $n$ is slightly large. It is indicated in Fig. 1 that the problem has the exponential ill-posedness.

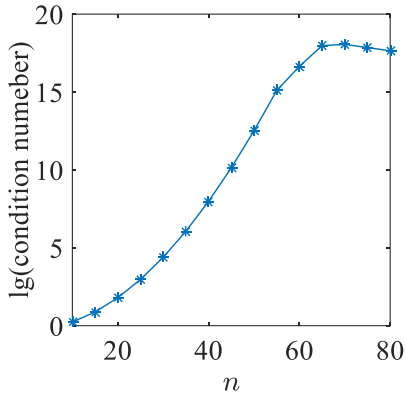

Fig.1 Dependence of the condition number of the matrix $\boldsymbol{K}$ on $n$

For having the data in the right side of Eq. (53), we assume the exact solution $\boldsymbol{f}^\dagger \equiv 1$. Then $\boldsymbol{d} = \boldsymbol{Kf}^\dagger$, and the noisy data is constructed through

$$[\boldsymbol{d}^\delta]_i = (1 + 2*(\text{rand}(1) - 0.5)*\delta')[\boldsymbol{d}]_i$$
$$i = 1, \cdots, n$$

Then the noise level of the measurement data is calculated by $\delta = \|\boldsymbol{d}^\delta - \boldsymbol{d}\|$. With matrix $\boldsymbol{K}$ and data $\boldsymbol{d}^\dagger$, the iterative methods (41), (43)-(44), and (47) are applied to solve the noisy system

$$\boldsymbol{Kf} = \boldsymbol{d}^\delta \quad (54)$$

In all methods, we simply set the initial values $\boldsymbol{f}_0 = 0$, $\dot{\boldsymbol{f}}_0 = 0$ and the magnitude $n = 100$. For an approximate solution $\boldsymbol{f}^k$, its accuracy is assessed by using the relative error as

$$\text{L2err} = \|\boldsymbol{f}^k - \boldsymbol{f}^\dagger\| / \|\boldsymbol{f}^\dagger\|$$

In addition, IterN is used to denote iterative number where the iteration stops.

The effect of the step-size $\Delta t$ on the iterative number and the solution accuracy is firstly investigated. To the end, fix $\delta' = 1\%, \tau = 1.03$. The results are shown in Figs. 2 and 3. Fig. 2 corresponds to constant damping parameters ($\eta = 0.6$ for SE1, $\eta = 0.8$ for SV1, $\eta = 0.1$ for MSV1 and $\eta = 0.1$ for RK1) while Fig. 3 corresponds to a variable damping parameter ($\eta(t) = 4/t$ for all methods). It is concluded that for all methods, on one hand, the larger $\Delta t$ is, the faster the iteration is; on the other hand, its too big values may lead to the divergence. In the following experiments, we set $\Delta t = 0.7$ for SE1, $\Delta t = 0.8$ for SV1, $\Delta t = 0.4$ for MSV1, $\Delta t = 1.1$ for RK1, $\Delta t = 0.6$ for SE2, $\Delta t = 0.8$ for SV2, $\Delta t = 0.4$ for MSV2, $\Delta t = 1.1$ for RK2. They are all approximately optimal.

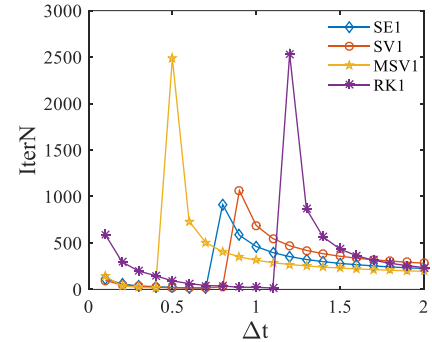

(a) Dependence of IterN on $\Delta t$

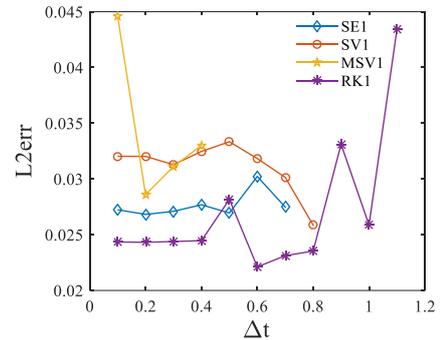

(b) Dependence of L2err on $\Delta t$

Fig.2 Dependence of IterN and L2err on $\Delta t$ with $\eta$ as a constant

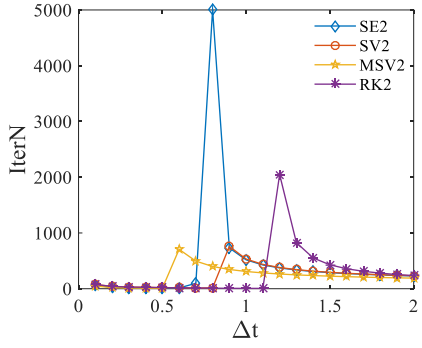

(a) Dependence of IterN on $\Delta t$

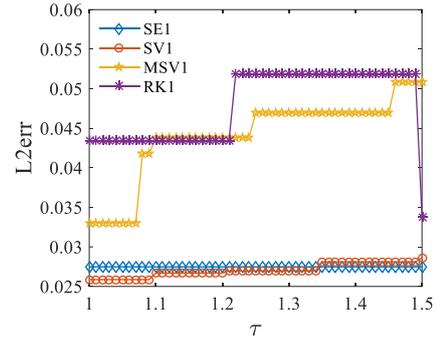

(b) Dependence of L2err on $\tau$

Fig.4 Dependence of IterN and L2err on $\tau$ with $\eta$ as a constant

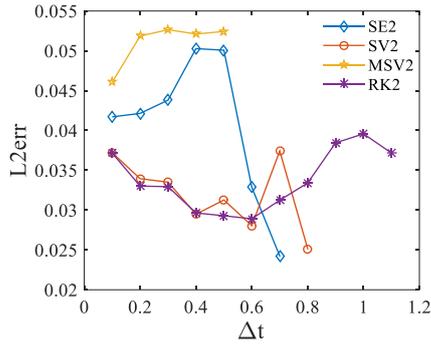

(b) Dependence of L2err on $\Delta t$

Fig.3 Dependence of IterN and L2err on $\Delta t$ with $\eta(t) = 4/t$

Next we investigate the effect of the parameter $\tau$, used in discrepancy principle, on the iterative number and the solution accuracy. For this purpose, fix $\delta' = 1\%$. The results are plotted in Figs. 4 and 5 which show that on the whole, the larger the value of $\tau$ is, the less the iterative number is and the worse the solution accuracy is. In the remaining part, without specific statement, we set $\tau = 1.03$.

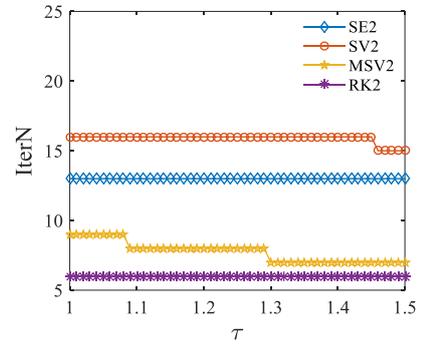

(a) Dependence of IterN on $\tau$

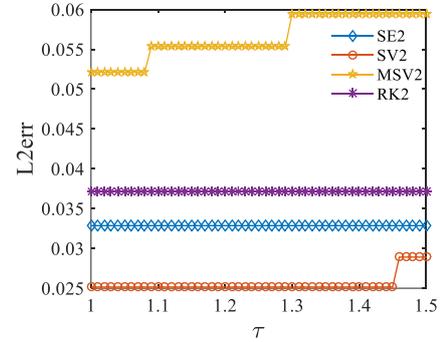

(b) Dependence of L2err on $\tau$

Fig.5 Dependence of IterN and L2err on $\tau$ with $\eta(t) = 4/t$

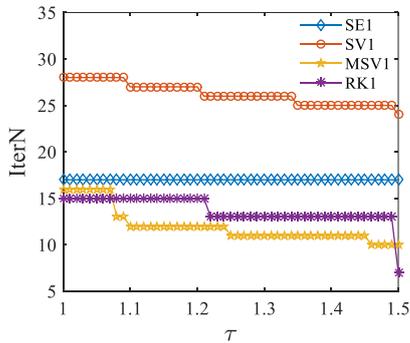

(a) Dependence of IterN on $\tau$

We further compare the iterative schemes with the existing Landweber method, CG method, $\nu$-method and Nesterov method. We set $\Delta t = 0.3$ for Landweber method, $\alpha = 3, \omega = 0.16$ for Nesterov method, $\Delta t = 0.7$, $\eta = 0.6$ for SE1, $\Delta t = 0.8, \eta = 0.8$ for SV1, $\Delta t = 0.4, \eta = 0.1$ for MSV1, $\Delta t = 1.1, \eta = 0.1$ for RK1, $\Delta t = 0.6$ for SE2, $\Delta t = 0.8$ for SV2, $\Delta t = 0.4$ for MSV2, $\Delta t = 1.1$ for RK2. Moreover, set $\tau = 1.03$ in all methods. The iterative numbers and relative errors in approximate solutions for different methods and different noise levels are given in

Table 1, from which we conclude that with properly chosen parameters, all the mentioned methods are stable and can produce satisfactory solutions. Compared with the Landweber method, all the other mentioned methods have acceleration effect. On the whole, for all methods, the larger the noise level is, the worse the solution accuracy is, but the less the required iterative number is. Compared with SE1, SV1, MSV1 and RK1 methods, SE2, SV2, MSV2 and RK2 have better behavior, such as fewer iterative steps or better solution accuracy. Particularly, Table 1 shows that RK2 possesses the best behavior in both the iterative number and solution accuracy.

**Table 1 Comparison between different methods in Example 1**

| $\delta'$ | 0.1% | | 1% | | 5% | |
|---|---|---|---|---|---|---|
| Example 1 | | | | | | |
| Methods | L2err | IterN | L2err | IterN | L2err | IterN |
| Landweber | 2.2102e−2 | 112 | 3.5974e−2 | 28 | 6.6438e−2 | 19 |
| CG | 1.8990e−2 | 7 | 4.4845e−2 | 3 | 9.7072e−2 | 2 |
| $v=0.5$ | 2.2606e−2 | 541 | 2.3479e−2 | 55 | 7.6998e−2 | 12 |
| $v=0.7$ | 9.6549e−3 | 99 | 2.6060e−2 | 19 | 7.3565e−2 | 6 |
| $v=1.0$ | 1.7758e−2 | 33 | 3.7817e−2 | 9 | 7.0558e−2 | 4 |
| $v=1.5$ | 2.2145e−2 | 25 | 4.9879e−2 | 6 | 7.3677e−2 | 3 |
| $v=2.0$ | 2.2240e−2 | 27 | 5.1933e−2 | 6 | 7.4185e−2 | 3 |
| Nesterov | 2.0208e−2 | 44 | 4.9961e−2 | 9 | 8.6451e−2 | 3 |
| SE1 | 2.1083e−2 | 34 | 2.7476e−2 | 17 | 6.7409e−2 | 16 |
| SV1 | 1.9937e−2 | 49 | 2.5863e−2 | 28 | 6.9638e−2 | 16 |
| MSV1 | 1.5788e−2 | 52 | 3.3020e−2 | 16 | 8.1759e−2 | 3 |
| RK1 | 1.0472e−2 | 49 | 4.3431e−2 | 15 | 8.0282e−2 | 7 |
| SE2 | 1.4591e−2 | 56 | 3.2859e−2 | 13 | 6.3635e−2 | 8 |
| SV2 | 1.3030e−2 | 53 | 2.5152e−2 | 16 | 6.9788e−2 | 14 |
| MSV2 | 2.2121e−2 | 40 | 5.2112e−2 | 9 | 8.7252e−2 | 4 |
| RK2 | 2.0732e−2 | 16 | 3.7170e−2 | 6 | 6.4510e−2 | 5 |

Finally we study the effect of the magnitude of the problem on different iterative schemes. We fix $\delta'=1\%$, $\tau=1.03$. The experiments for Landweber method, $v$-method, Nesterov method and RK2 are implemented repeatedly to solve Eq. (54) with $n=25, 50, 100, 200, 400, 800, 1600$ and $3200$, respectively, and the results are shown in Fig. 6. It is indicated from the Fig. 6 that when the order of the matrix increases, the accuracy of the solution and the number of iteration steps change little. Moreover, Fig. 6 also makes clear that the RK2 performs better than three existing methods, especially when the scale of the problem increases.

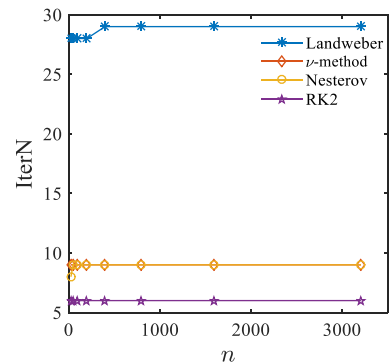

(a) Dependence of IterN on the scale $n$

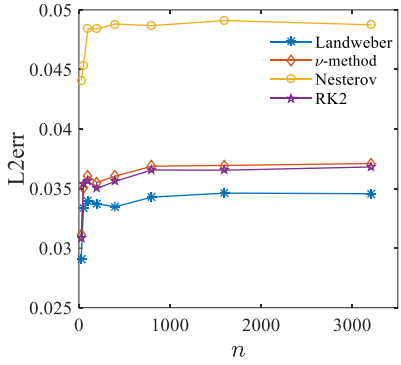

(b) Dependence of L2err on the scale $n$

Fig.6 Dependence of IterN and L2err on the scale $n$

### 3.2. Example 2

In the second example, the well-known ill-posed Hilbert matrix is taken to model the operator as

$$A = \begin{pmatrix} 1 & \frac{1}{2} & \frac{1}{3} & \cdots & \frac{1}{n} \\ \frac{1}{2} & \frac{1}{3} & \frac{1}{4} & \cdots & \frac{1}{n+1} \\ \cdots & \cdots & \cdots & \cdots & \cdots \\ \frac{1}{n} & \frac{1}{n+1} & \frac{1}{n+2} & \cdots & \frac{1}{2n-1} \end{pmatrix}$$

The dependence of the condition number of the matrix $A$ on its scale is plotted in Fig. 7, which shows that the condition number of $A$ increases dramatically and thus $A$ is ill-posed. Moreover, like Example 1, assume again the exact solution $x^\dagger = (1,1,\cdots,1)^T$ and compute $Ax^\dagger$ for the right side $b$. The noisy data $b^\delta$ is constructed like Example 1. Then iterations (41), (43)-(44) and (47) are applied to solve

$$Ax = b^\delta \tag{55}$$

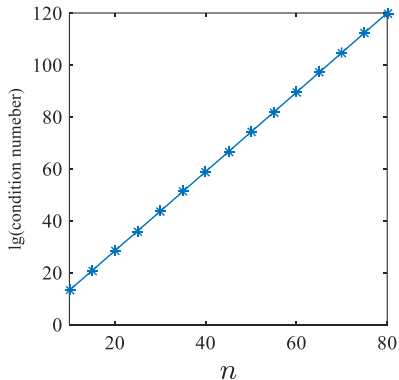

Fig.7 Dependence of condition number of $A$ on the scale $n$

The dependence of the iterative numbers and the accuracy in approximate solutions obtained with iterations in Section 3 on parameters $\Delta t$ and $\tau$ are plotted in Figs. 8-11. Specifically, Fig. 8 is for constant $\eta$ (with $\eta = 0.2$ for SE1, $\eta = 0.2$ for SV1, $\eta = 0.1$ for MSV1, $\eta = 0.1$ for RK1) while Fig. 9 is for $\eta(t) = 4/t$ (both with $\tau = 1.03$); Fig. 10 is for constant $\eta$ (with $\Delta t = 0.8$, $\eta = 0.2$ for SE1, $\Delta t = 0.9$, $\eta = 0.2$ for SV1, $\Delta t = 0.5, \eta = 0.1$ for MSV1, $\Delta t = 1.2$, $\eta = 0.1$ for RK1) while Fig. 11 is for $\eta(t) = 4/t$ (with $\Delta t = 0.7$ for SE2, $\Delta t = 0.9$ for SV2, $\Delta t = 0.5$ for MSV2, $\Delta t = 1.1$ for RK2). In all figures, $\delta' = 1\%$. Similar conclusion to those for Example 1 can be drawn from these figures.

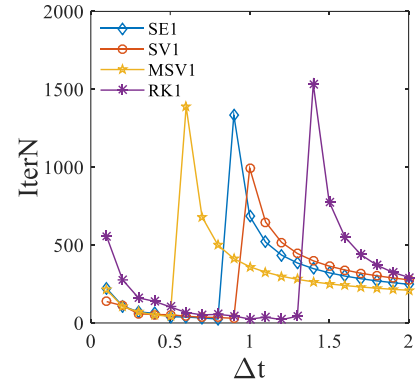

(a) Dependence of IterN on $\Delta t$

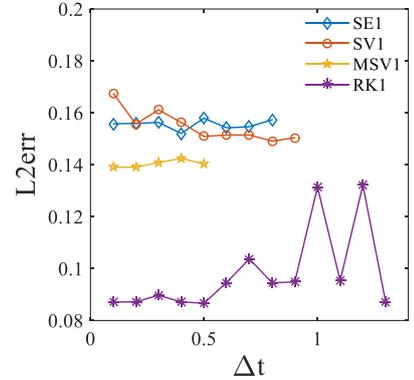

(b) Dependence of L2err on $\Delta t$

Fig.8 Dependence of IterN and L2err on $\Delta t$ with $\eta$ as a constant

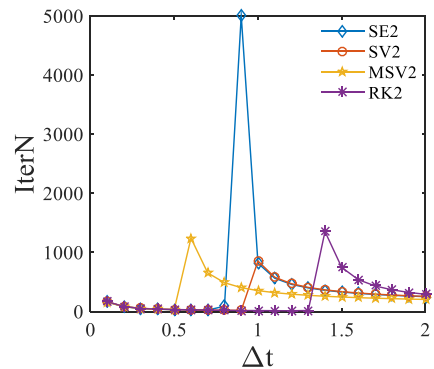

(a) Dependence of IterN on $\Delta t$

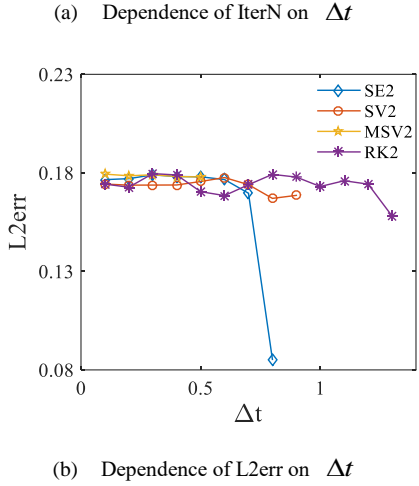

(b) Dependence of L2err on $\Delta t$

Fig.9 Dependence of IterN and L2err on $\Delta t$ with $\eta(t)=4/t$

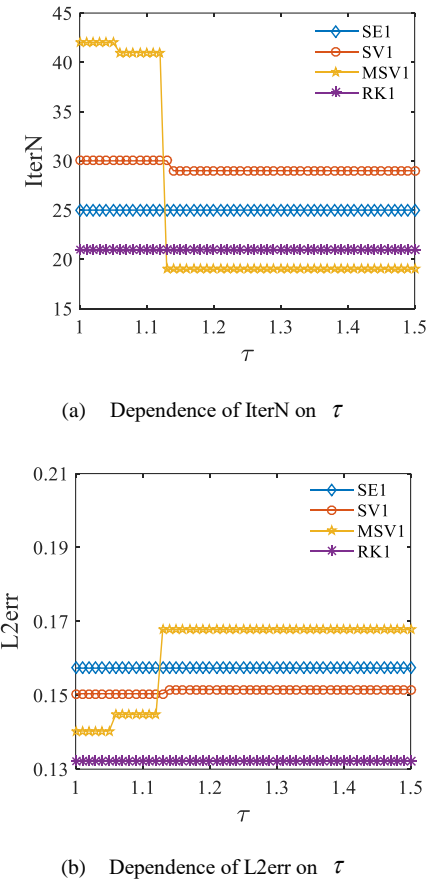

(a) Dependence of IterN on $\tau$

(b) Dependence of L2err on $\tau$

Fig.10 Dependence of IterN and L2err on $\tau$ with $\eta$ as a constant

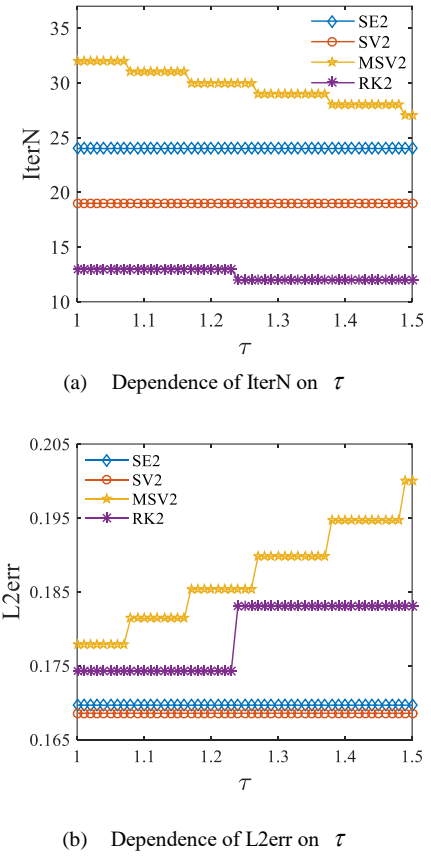

(a) Dependence of IterN on $\tau$

(b) Dependence of L2err on $\tau$

Fig.11 Dependence of IterN and L2err on $\tau$ with $\eta(t)=4/t$

Similar to those in Example 1, a comparison between different methods for different noise data is displayed in Table 2, where $\Delta t=0.3$ for Landweber method, $\alpha=3$, $\omega=0.2$ for Nesterov method, $\Delta t=0.8, \eta=0.2$ for SE1, $\Delta t=0.9, \eta=0.2$ for SV1, $\Delta t=0.5, \eta=0.1$ for MSV1, $\Delta t=1.2, \eta=0.1$ for RK1, $\Delta t=0.7$ for SE2, $\Delta t=0.9$ for SV2, $\Delta t=0.5$ for MSV2, $\Delta t=1.1$ for RK2. Also, $\tau=1.03$ in all methods. It can be concluded from Table 2 that with properly chosen parameters, all the mentioned methods are stable and can produce satisfactory solutions.

**Table 2 Comparison between different methods in Example 2**

| $\delta'$ | 0.1% | | 1% | | 5% | |
|---|---|---|---|---|---|---|
| Example 2 | | | | | | |
| Methods | L2err | IterN | L2err | IterN | L2err | IterN |
| Landweber | 8.1242e−2 | 2461 | 1.7665e−1 | 126 | 3.4507e−1 | 11 |
| CG | 6.9416e−2 | 4 | 1.3891e−1 | 3 | 2.9184e−1 | 2 |
| $v=0.5$ | 2.9539e−2 | 545 | 8.4796e−2 | 60 | 2.0654e−1 | 12 |
| $v=0.7$ | $7.4753\text{e}-2$ | 130 | $1.4767\text{e}-1$ | 28 | $2.9460\text{e}-1$ | 8 |
| $v=1.0$ | $8.1044\text{e}-2$ | 82 | $1.7580\text{e}-1$ | 18 | $3.4722\text{e}-1$ | 5 |
| $v=1.5$ | $8.1710\text{e}-2$ | 98 | $1.7892\text{e}-1$ | 21 | $3.4342\text{e}-1$ | 6 |
| $v=2.0$ | $8.1712\text{e}-2$ | 113 | $1.7500\text{e}-1$ | 25 | $3.4880\text{e}-1$ | 6 |
| Nesterov | $7.8666\text{e}-2$ | 159 | $1.7555\text{e}-1$ | 32 | $3.4949\text{e}-1$ | 8 |
| SE1 | $8.1388\text{e}-2$ | 180 | $1.4402\text{e}-1$ | 29 | $1.2750\text{e}-1$ | 29 |
| SV1 | $8.1414\text{e}-2$ | 160 | $1.2063\text{e}-1$ | 46 | $1.3307\text{e}-1$ | 28 |
| MSV1 | $8.1903\text{e}-2$ | 136 | $1.4033\text{e}-1$ | 42 | $2.6692\text{e}-1$ | 11 |
| RK1 | $7.9167\text{e}-2$ | 66 | $1.3216\text{e}-1$ | 21 | $2.3541\text{e}-1$ | 5 |
| SE2 | $7.8706\text{e}-2$ | 114 | $1.6966\text{e}-1$ | 24 | $2.5410\text{e}-1$ | 14 |
| SV2 | $7.9494\text{e}-2$ | 86 | $1.6857\text{e}-1$ | 19 | $2.7046\text{e}-1$ | 10 |
| MSV2 | $8.1590\text{e}-2$ | 148 | $1.7791\text{e}-1$ | 32 | $3.5131\text{e}-1$ | 8 |
| RK2 | $8.1213\text{e}-2$ | 67 | $1.7600\text{e}-1$ | 14 | $3.1168\text{e}-1$ | 6 |

The effect of the problem magnitude on different iterative schemes is also studied. The experiments for Landweber method, $v$-method, Nesterov method and RK2 are implemented repeatedly to solve Eq. (55) with $\delta'=1\%, \tau=1.03$, and $n=25, 50, 100, 200, 400, 800, 1600$ and 3200, respectively, whose results are shown in Fig. 12. Again we can see from Fig. 12 that the RK2 performs the best when compared with the existing classical methods.

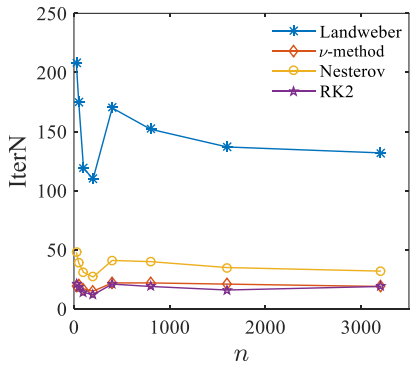

(a) Dependence of IterN on the scale $n$

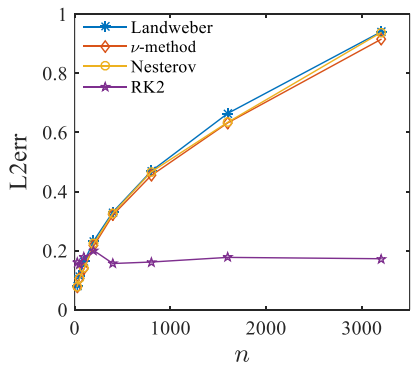

(b) Dependence of L2err on the scale $n$

Fig.12 Dependence of IterN and L2err on the scale $n$

## 4  Conclusion

The newly developed second order asymptotic system is applied to solve severely ill-posed linear system. Compared with the existing work, a series of theoretical results are presented and proved with simplified arguments. Since the existing work for second order system often assumes that it has exact solution while our problem may have no solution, the analysis is adjusted. On the whole, the fourth order Runge-Kutta method is better than the classical $v$-method and Nesterov method, and much faster than Landweber method. Moreover, it is more obvious when the scale and the ill-posedness of the system increase. The reason is that the system is ill-posed and iteration should stop in advance before the accuracy gets worse. Therefore, although it is not a symplectic method, it has high precision. This is a new sight. Better behavior can be expected if higher precision scheme is used. The idea can also be applied to other linear inverse problems such as inverse source problems and Cauchy problems etc., which will be studied in the future.